\documentclass[letterpaper, 10 pt, conference]{ieeeconf}  

\IEEEoverridecommandlockouts                              

\usepackage{graphics} 
\usepackage{epsfig} 
\usepackage{color}
\usepackage{mathptmx} 
\usepackage{times} 
\usepackage{amsmath} 
\usepackage{amssymb}  
\usepackage{todonotes}
\usepackage{graphicx}
\usepackage{mathtools}
\usepackage{dirtytalk}
\usepackage[caption=false]{subfig}

\mathchardef\mhyphen="2D 

\title{\LARGE \bf
The Parallelization of Riccati Recursion 
}

\author{Forrest Laine$^{1}$ and Claire Tomlin$^{1}$
\thanks{This research is supported by DARPA under the Assured Autonomy Program, 
by NSF under the CPS Frontier VeHICal project, and by the UC-Philippine-California Advanced Research Institute under project IIID-2016-005.}
\thanks{$^{1}$Forrest Laine and Claire Tomlin are with the Department of Electrical Engineering and Computer Sciences, 
        University of California, Berkeley
        {\tt\small forrest.laine@berkeley.edu}, \tt\small tomlin@eecs.berkeley.edu}%
        }

\begin{document}

\maketitle
\thispagestyle{empty}
\pagestyle{empty}

\begin{abstract}
A method is presented for parallelizing the computation of solutions to discrete-time, linear-quadratic, finite-horizon optimal control problems, which we will refer to as LQR problems. This class of problem arises frequently in robotic trajectory optimization. For very complicated robots, the size of these resulting problems can be large enough that computing the solution is prohibitively slow when using a single processor. Fortunately, approaches to solving these type of problems based on numerical solutions to the KKT conditions of optimality offer a parallel solution method and can leverage multiple processors to compute solutions faster. However, these methods do not produce the useful feedback control policies that are generated as a by-product of the dynamic-programming solution method known as Riccati recursion. In this paper we derive a method which is able to parallelize the computation of Riccati recursion, allowing for super-fast solutions to the LQR problem while still generating feedback control policies. We demonstrate empirically that our method is faster than existing parallel methods. 


\end{abstract}

\section{INTRODUCTION} \label{sec:intro}
When solving a continuous-state, continuous-action trajectory optimization problem for robotic path planning, many of the well known methods for solving these problems result in subproblems taking the form of a discrete-time, time-varying, finite-horizon optimal control problem with linear dynamics and a quadratic objective \cite{bertsekas1999nonlinear}, \cite{wright1993interior}, \cite{wright1996applying}, \cite{wright1999numerical} \cite{allgower2012nonlinear}. The class of such problems has a long history of research, and is commonly referred to as the Linear Quadratic Regulator (LQR) \cite{kwakernaak1972linear}.  In applications in which the trajectory optimization must be computed online, such as in Model Predictive Control frameworks \cite{allgower2012nonlinear}, we require methods for solving this subproblem very efficiently.  

Fortunately, LQR problems offer efficient methods for computing solutions based on dynamic programming, sometimes referred to as the discrete-time Riccati recursion. Alternatively, one could solve these types of problem by direct methods based on the KKT system of equations given for the optimization problem. Both of these methods require computation that scales linearly with respect to the number of discrete time-points in the trajectory \cite{borrelli2017predictive}, \cite{wright1996applying}. The direct methods based on solving the KKT system of equations rely on the banded nature of the coefficient matrix of the equations. For this reason, the solution to these problems can also utilize parallelized methods for solving banded systems of equations \cite{dongarra1987solving}, \cite{arbenz1999comparison}.  These parallel methods allow applications to utilize all of the available cores available on a computer to compute solutions as fast as possible. 

Though there exist parallelized methods for computing optimal trajectories in the LQR problem, because they are based on direct methods, they only compute numerical values of the trajectories and do not include the feedback control policies that are generated through Riccati recursion. In online Model Predictive Control applications, these policies can be used to stabilize the system to the planned trajectory until the next trajectory is generated. These policies can also be used to stabilize single-shooting approaches used in nonlinear optimization methods \cite{giftthaler2017family}. 

It would therefore be desirable to have a method which can parallelize the solution to Riccati recursion, so that feedback control policies can be generated in a way that also leverages multiple processors to divide up computation. There are two main ingredients that are needed to derive such a method. First, a means of breaking up the trajectory optimization problem, which we will call the \textit{global} problem, into subproblems so that they can be computed in parallel. Second, a method is needed which can solve those subproblems such that their solutions can be reconstructed into the solution of the global problem.  The contributions of this paper are derivations of both ingredients, and as such, a resulting method for parallelizing Riccati recursion.

In section \ref{sec:method} we present our method. In that presentation, we will first introduce the global problem, and propose a method for how it can be parallelized. We then will present a means for computing the solution to the resulting subproblems, and show how the solutions can be reconstructed. In section \ref{sec:compare} we discuss properties of our method, and demonstrate empirically that our method in fact computes the solutions to the global problem. We also show that our method is faster than the fastest-known method for solving the KKT system of equations in parallel. 

\section{PARALLELIZING RICCATI RECURSION} \label{sec:method}

The global LQR problem which we wish to parallelize will be defined as the following:
\begin{subequations}
\begin{align}
    & \min_{x_0, u_0,...,u_{T-1}, x_T} cost_T(x_T) + \sum_{t=0}^{T-1} cost_t(x_t, u_t) \label{obj:globalobjwords} \\
    & \text{s.t.} \ \ \ dynamics_t(x_{t+1}, x_t, u_t) = 0 \ \forall t \in \{0,...,T-1\} \label{eq:dynamicswords} \\
    & \ \ \ \ \  \ \ x_0 = x_{init} \label{eq:init} 
\end{align} \label{opt:globalwords}
\end{subequations}
Where $x_t \in \mathbb{R}^n$, $u_t \in \mathbb{R}^m$,  the functions $cost_t: \mathbb{R}^n \times \mathbb{R}^m \to \mathbb{R}$, $cost_T: \mathbb{R}^n \to \mathbb{R}$, $dynamics_t : \mathbb{R}^n \times \mathbb{R}^n \times \mathbb{R}^m \to \mathbb{R}^n$, are functions defined as:
\begin{subequations}
\begin{align}
    &cost_t(x, u) = \frac{1}{2} \begin{pmatrix} 1 \\ x \\ u \end{pmatrix}^\intercal \begin{pmatrix} 0 & q_{x1_t}^\intercal & q_{u1_t}^\intercal \\ q_{x1_t} & Q_{xx_t} & Q_{ux_t}^\intercal \\ q_{u1_t} & Q_{ux_t} & Q_{uu_t} \end{pmatrix}\begin{pmatrix} 1 \\ x \\ u \end{pmatrix} \label{eq:costtform} \\
    &cost_T(x) = \frac{1}{2} \begin{pmatrix} 1 \\ x \end{pmatrix}^\intercal \begin{pmatrix} 0 & q_{x1_T}^\intercal \\ q_{x1_T} & Q_{xx_T} \end{pmatrix}\begin{pmatrix} 1 \\ x \end{pmatrix} \\
    &dynamics_t(x_{t+1}, x_t, u_t) = x_{t+1} - (F_{x_t} x_t + F_{u_t} u_t + f_{1_t}) \label{eq:dynamicsform}
\end{align}
\end{subequations}

We assume that the block coefficients $Q_{uu_t}$ of the quadratic functions $cost_t$ are positive-definite, and the terms $ Q_{xx_t}-Q_{ux_t}Q_{uu_t}^{-1}Q_{ux_t}^\intercal$ is positive semi-definite, so that problem (\ref{opt:globalwords}) is strictly convex. 

\subsection{Dividing into Sub-trajectories}

\begin{figure}[t] 
\centering
\subfloat[Step 1: Compute sub-trajectories and Lagrange multipliers as functions of unknown link-points.]{
  \centering
  \includegraphics[scale=0.12]{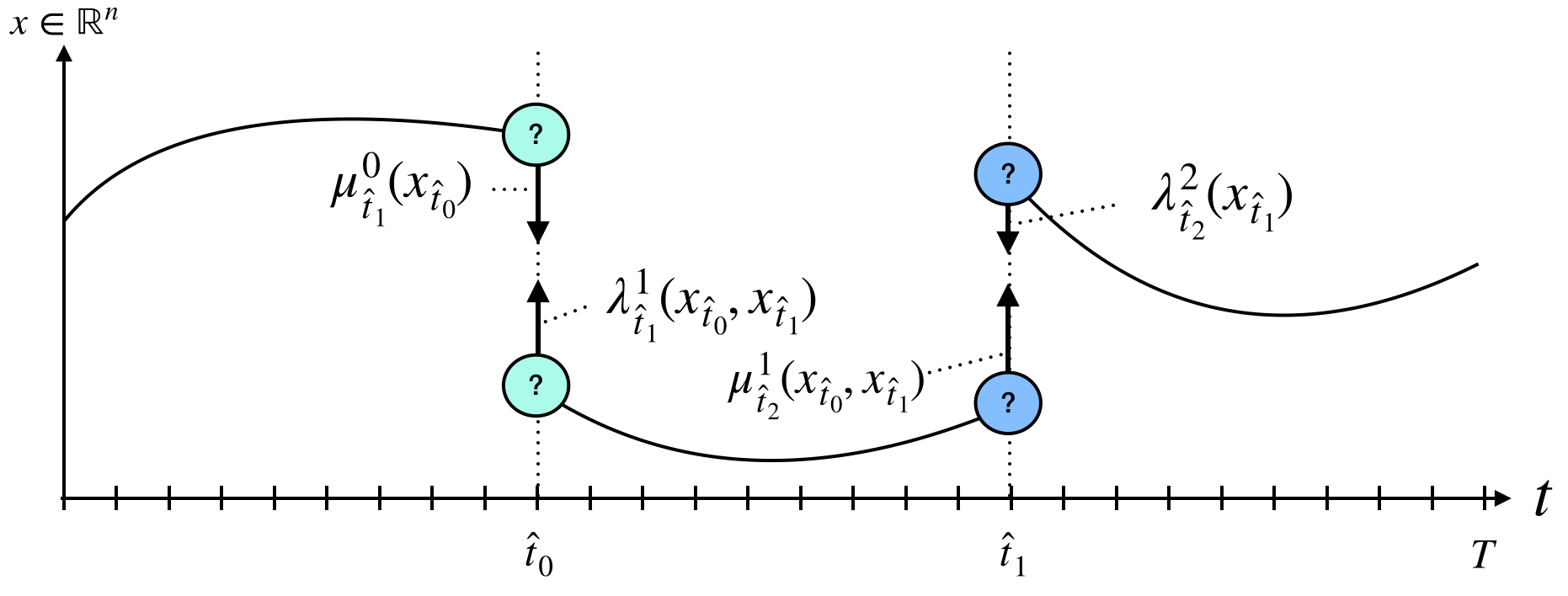} \label{fig:pre_combine}
  }

\subfloat[Step 2: Relate neighboring Lagrange multipliers to solve for link-points.] {
  \centering
  \includegraphics[scale=0.12]{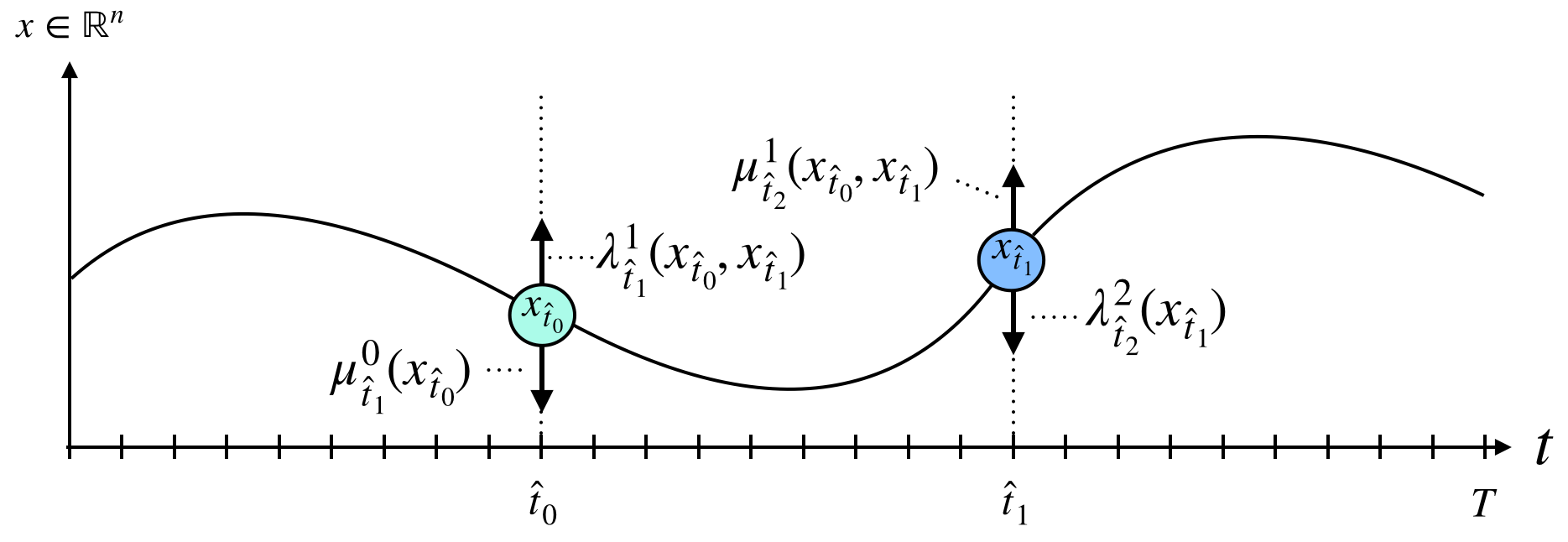} \label{fig:post_combine}
  }
  \caption{A depiction of the parallelization of an LQR computation.} 
  \label{fig:cartoon}
\end{figure}

The way problem (\ref{opt:globalwords}) will be broken into independent sub-problems will come from the KKT conditions of optimality. Together, along with the dynamic constraints (\ref{eq:dynamicswords}), the KKT conditions for problem (\ref{opt:globalwords}) can be stated as 

\begin{subequations}
\begin{align}
&\begin{aligned}
\begin{bmatrix} Q_{xx_t} & Q_{ux_t}^\intercal \\  Q_{ux_t} & Q_{uu_t} \end{bmatrix} \begin{bmatrix} x_t \\ u_t \end{bmatrix} + \begin{bmatrix} q_{x1_t} \\ q_{u1_t} \end{bmatrix} + \begin{bmatrix} I & -F_{x_t}^\intercal \\ & -F_{u_t}^\intercal \end{bmatrix} \begin{bmatrix} \lambda_t \\ \lambda_{t+1} \end{bmatrix} = \begin{bmatrix} 0 \\ 0 \end{bmatrix} \\
 \ \ \ \ \ \ \ \ \  \forall t \in \{0,...,T-1\} 
\end{aligned} \label{eq:mults1} \\
&Q_{xx_T} x_T^* + q_{x1_T} + \lambda_{T}  = 0 \ \ \ \ \ \ \ \ \ \ \ \ \ \ \ \ \ \ \ \ \ \ \ \ \ \label{eq:multsT1}
\end{align} \label{eq:KKT}
\end{subequations}

Here $\lambda_{t+1}$ is the Lagrange multiplier associated with constraints (\ref{eq:dynamicswords}), and $\lambda_0$ is the multiplier associated with constraint (\ref{eq:init}). Finding a set of $x_t$, $u_t$, $\lambda_t$ that satisfy these equations will result in the optimal solution to the global problem \cite{boyd2004convex}. 

Notice that we can modify the equations (\ref{eq:KKT}) around some time $\tau$ to obtain another set of equivalent equations: 

\begin{subequations}
 \begin{align}
&\begin{aligned}
\begin{bmatrix} Q_{xx_t} & Q_{ux_t}^\intercal \\  Q_{ux_t} & Q_{uu_t} \end{bmatrix} \begin{bmatrix} x_t \\ u_t \end{bmatrix} + \begin{bmatrix} q_{x1_t} \\ q_{u1_t} \end{bmatrix} + \begin{bmatrix} I & -F_{x_t}^\intercal \\ & -F_{u_t}^\intercal \end{bmatrix} \begin{bmatrix} \lambda_t \\ \lambda_{t+1} \end{bmatrix} = \begin{bmatrix} 0 \\ 0 \end{bmatrix} \\
 \ \ \ \ \ \ \ \ \  \forall t \in \{0,...,\tau-1\} 
\end{aligned} \label{eq:mults1} \\
& \lambda_{\tau} + \mu_{\tau} = 0 \ \ \ \ \ \ \ \ \ \ \ \ \ \ \ \ \ \ \ \ \ \ \ \ \ \label{eq:mu1} \\
& x_{\tau} = x_{link} \ \ \ \ \ \ \ \ \ \ \ \ \ \ \ \ \ \ \ \ \ \ \ \ \\
&\begin{aligned}
\begin{bmatrix} Q_{xx_t} & Q_{ux_t}^\intercal \\  Q_{ux_t} & Q_{uu_t} \end{bmatrix} \begin{bmatrix} x_t \\ u_t \end{bmatrix} + \begin{bmatrix} q_{x1_t} \\ q_{u1_t} \end{bmatrix} + \begin{bmatrix} I & -F_{x_t}^\intercal \\ & -F_{u_t}^\intercal \end{bmatrix} \begin{bmatrix} \lambda_t \\ \lambda_{t+1} \end{bmatrix} = \begin{bmatrix} 0 \\ 0 \end{bmatrix} \\
 \ \ \ \ \ \ \ \ \  \forall t \in \{\tau,...,T-1\} 
\end{aligned} \label{eq:mults2} \\
&Q_{xx_T} x_T+ q_{x1_T} + \lambda_{T} = 0 \ \ \ \ \ \ \ \ \ \ \ \ \ \ \ \ \ \ \ \ \ \ \ \ \ \label{eq:mu2}
\end{align} \label{eq:KKT2}
\end{subequations}

Here, $x_{link} = x_\tau^*$, is the optimal point from the global problem. This has the effect of constraining $x_\tau$ to be its optimal value. The term $\mu_\tau$ is an auxiliary variable whose purpose will be come clear shortly. Notice that these set of equations (\ref{eq:KKT2}) correspond to the KKT conditions for two independent optimization problems:

\begin{subequations}
\begin{align}
    & \min_{x_0, u_0,...,u_{\tau-1}, x_\tau^1} \ \ \  \sum_{t=0}^{\tau-1} cost_t(x_t, u_t)  \label{eq:opt1} \\
    & \text{s.t.} \ \ \ dynamics_t(x_{t+1}, x_t, u_t) = 0 \ \forall t \in \{0,...,\tau-1\} \\
    & \ \ \ \ \  \ \ x_0 = x_{init} \label{eq:init1}  \\
    & \ \ \ \ \  \ \ x_\tau^1 = x_{link} \label{eq:term1} \\
    & \min_{x_\tau^2, u_\tau,...,u_{T-1}, x_T} cost_T(x_T) + \sum_{t=\tau}^{T-1} cost_t(x_t, u_t)  \\
    & \text{s.t.} \ \ \ dynamics_t(x_{t+1}, x_t, u_t) = 0 \ \forall t \in \{\tau,...,T-1\}  \\
    & \ \ \ \ \  \ \ x_\tau^2 = x_{link} \label{eq:init2} 
\end{align} \label{opt:split}
\end{subequations}

$x_\tau^1$ and $x_\tau^2$ are used to distinguish the variables in the two optimization problems, since $x_\tau$ appears in both problems. From this, it can be seen that the term $\mu_\tau$ appearing in equation (\ref{eq:mu1}) is the Lagrange multiplier corresponding to constraint (\ref{eq:term1}). 

Therefore, if the value of $x_{link}$ could be determined somehow, solving the problems (\ref{opt:split}) independently would provide a means for solving the global problem in parallel. It turns out that by computing the optimal solutions \textit{and} the Lagrange multipliers corresponding to the endpoint constraints in each of the sub-problems (\ref{opt:split}) as functions of the unknown point $x_{link}$, the resulting solutions can be used to solve for the point $x_{link}$. Since the solutions of the sub-problems are functions of $x_{link}$, by plugging it back into those solutions, the global solution can be recovered. 

The method for computing solutions as a function of the endpoint will presented in the following section. For now assume that such a method exists, giving solutions to each of the subproblems (\ref{opt:split}) as:

\begin{subequations}
\begin{align}
	x_\tau^1 &= x_\tau^2= x_{link} \\
	 \lambda_\tau^1 &= L_\tau^1 x_{link} + l_\tau^1 \\
	 \lambda_\tau^2 &= L_\tau^2 x_{link} + l_\tau^2 \\
	x_t^* &= R_t x_{link} + r_t, \ \ t \in \{0,...,\tau-1,\tau+1,...,T\} \\
	u_t^* &= S_t x_{link} + s_t, \ \ t \in \{0,...,T-1\} \\
	 \lambda_t &= L_t x_{link} + l_t, \ \ t \in \{0,...,\tau-1,\tau+1,...,T\} \\
	\mu_\tau &= E_\tau x_{link} + e_{\tau}.
\end{align} \label{eq:explicits}
\end{subequations}

By noticing that the condition  $\mu_\tau=-\lambda_\tau^2$ is required if the KKT conditions are to be equivalent to the equations in (\ref{eq:KKT}, \ref{eq:KKT2}), we obtain:
\begin{equation}
	(E_\tau + L_\tau^2) x_{link} + (e_\tau + l_\tau^2) = 0
\end{equation}
This is a system of $n$ equations in $n$ unknowns, which can be solved for $x_{link}$. This assumes $E_\tau + L_\tau^2$ is full rank, which it will be when the global problem is strictly convex. With the solution of $x_{link}$, the solutions of the subproblems (and by virtue of the equivalence, the global problem) given in (\ref{eq:explicits}) can be evaluated. 

This argument made relating the Lagrange multipliers of neighboring trajectories holds for an arbitrary number of sub-trajectories, given that a solution to each sub-optimization is feasible. When splitting the global problem into $J$ sub-problems, $J-1$ of those problems will have start and end-point constraints and no terminal cost, as in (\ref{eq:opt1}-\ref{eq:term1}), and one subproblem with a terminal cost and no endpoint constraint. Solutions of each sub-trajectory must be computed as function of all unknown link points, i.e. the values appearing in equations (\ref{eq:init1}-\ref{eq:term1}), if they are unknown. This will be the case for the $J-2$ sub-trajectories making up the middle of the global trajectory. 

The equivalence relationships needed to ensure the KKT conditions of all sub-problems match those of the global problem will result in $(J-1)n$ equations, with $(J-1)n$ unknowns, being the $J-1$ link-points. The system will be block-banded, with bandwidth equal to $2n+1$. Therefore the computational complexity of solving this combined system is $O(Jn^3)$.  Since typically $J<<T$, this cost of determining the link-points will be very small. 

Assuming there is an efficient method for computing the solutions (\ref{eq:explicits}) as explicit functions of the end-points, the overall method for computing solutions in the manner outlined will be very efficient. We present such a method now.

\subsection{Constraint-explicit LQR}

We present a method for computing control policies of the following end-point constrained LQR problem, explicitly as a function of the values $x_{init}$ and $x_{term}$.  

\begin{subequations}
\begin{align}
    & \min_{x_0, u_0,...,u_{T-1}, x_T} cost_T(x_T) + \sum_{t=0}^{T-1} cost_t(x_t, u_t) \label{obj:globalobjwords3} \\
    & \text{s.t.} \ \ \ dynamics_t(x_{t+1}, x_t, u_t) = 0 \ \forall t \in \{0,...,T-1\} \label{eq:dynamicswords3} \\
    & \ \ \ \ \  \ \ x_0 = x_{init} \label{eq:init3} \\ 
    & \ \ \ \ \  \ \ x_T = x_{term} \label{eq:term3}
\end{align} \label{opt:globalwords3}
\end{subequations}
Here the functions appearing in (\ref{obj:globalobjwords3}-\ref{eq:dynamicswords3}) take the same form as in (\ref{opt:globalwords}). It is assumed that the constraints (\ref{eq:dynamicswords3}-\ref{eq:term3}) are linearly independent. If there are linearly dependent constraints, they can simply be removed without changing the problem. 

The method we will use to solve this problem is based on the method presented in \cite{laine2018efficient}, except we will maintain all control policies, value-functions and constraints as explicit functions of the terminal state $x_{term}$. 

The derivation of the control policies for the constrained follows a dynamic programming approach. Starting from the end of the trajectory and working towards the beginning, a given control input $u_t$ will be chosen such that for any value of the corresponding state $x_t$ and terminal constraint value $x_{term}$, the remaining trajectory will be cost-optimal and the residual parts of the terminal constraint will be satisfied if possible.

To formalize and calculate these controls, a time-varying quadratic value function, $cost\_to\_go_t: \mathbb{R}^n \times \mathbb{R}^n \to \mathbb{R}$ will be maintained, representing the minimum possible cost remaining in the trajectory from stage $t$ onward as a function of state at time $t$ and the terminal constraint state $x_{term}$. Additionally, a linear function $constraint\_to\_go_t : \mathbb{R}^n \times \mathbb{R}^n \to \mathbb{R}^{r_t}$ will be maintained which encompasses all future constraints. Let the following terms represent these functions:
\begin{align}
\begin{aligned}
    cost\_to\_go_t(x, x_{term}) &= \\
    \ \ \ \ \ \ \ \  \frac{1}{2}
     \begin{pmatrix} 1 \\ x \\ x_{term} \end{pmatrix}^\intercal & \begin{pmatrix} 0 & v_{x1_t}^\intercal & v_{z1_t}^\intercal \\ v_{x1_t} & V_{xx_t} & V_{zx_t}^\intercal \\ v_{z1_t} & V_{zx_t} & V_{zz_t} \end{pmatrix}\begin{pmatrix} 1 \\ x \\ x_{term} \end{pmatrix} \label{eq:costtogo} 
\end{aligned} \\
constraint\_to\_go_t(x, x_{term}) = H_{x_t} x + H_{z_t}x_{term} + h_{1_t} .  \label{eq:constrainttogo}
\end{align}

We initialize these terms at time $T$ as follows: 
\begin{equation}
\begin{aligned}
    H_{x_T} &= I   & H_{z_T} &= -I \\
    h_{1_T} &= 0  & V_{xx_T} &= Q_{xx_T}  \\
    v_{x1_T} &= q_{x1_T}  & v_{z1_t} &= 0 \\
    V_{zx_T} &= 0 & V_{zz_T} &= 0
\end{aligned}
\end{equation}
Here we abuse notation slightly, letting $0$ represent both the zero vector and matrix, with corresponding size which should be clear from context. $I$ is the $n\times n$ identity matrix. 

According to the procedure outlined in \cite{laine2018efficient}, the following optimization problem is solved starting at time $t=T-1$ and continuing back to $t=0$:

\begin{subequations} \begin{align}
    u_t^* = \ &\text{arg}\min_{u_t} \ cost_t(x_t, u_t) + cost\_to\_go_{t+1}(x_{t+1}, x_{term}) \\
    \text{s.t.} \ \ \ & 0 = dynamics_t(x_{t+1}, x_t, u_t)  \label{eq:subwordsdyn} \\
    & u_t \in  \text{arg}\min_u \| constraint\_to\_go_{t+1}(x_{t+1}, x_{term}) \|_2 \label{eq:leastsquareswords}
\end{align}  \label{opt:words} \
\end{subequations}
Again, as shown in \cite{laine2018efficient}, problem (\ref{opt:words}) is equivalent to the following problem:
\begin{subequations}
\begin{align}
&\begin{aligned}
y_t^*, w_t^* = \text{arg}\min_{v_t, w_t}  \frac{1}{2} \| N_{x_t} x_t + N_{u_t}P_{y_t} y_t + N_{z_t} x_{term} + n_{1_t} \|_2  \ + 
  \\  \frac{1}{2} \begin{pmatrix} 1 \\ x_t \\ Z_{w_t} w_t  \\ x_{term} \end{pmatrix}^\intercal  \begin{pmatrix} 0 & m_{x1_t}^\intercal  & m_{u1_t}^\intercal & m_{z1_t}^\intercal \\ m_{x1_t} & M_{xx_t} & M_{ux_t}^\intercal  & M_{zx_t}^\intercal \\ m_{u1_t} & M_{ux_t} & M_{uu_t} & M_{zu_t}^\intercal \\ m_{z1_t} & M_{zx_t} & M_{zu_t} & M_{zz_t} \end{pmatrix}\begin{pmatrix} 1 \\ x_t \\ Z_{w_t} w_t \\ x_{term} \end{pmatrix} \  \label{eq:vwsoln} \end{aligned} \\
&  \ \ \ \ \ u_t^* = P_{y_t} y_t^* + Z_{w_t} w_t^* \label{eq:udirectsum}
\end{align} \label{opt:infeasible}
\end{subequations}

Where the above terms are defined as:
\begin{align*}
    m_{x1_t} &= q_{x1_t} + F_{x_t}^\intercal v_{x1_{t+1}} & m_{u1_t} &= q_{u1_t} + F_{u_t}^\intercal v_{x1_{t+1}} \\ 
    m_{z1_t} &= v_{z1_t} & M_{zx_t} &= V_{zx_t} \\
    M_{zu_t} &= 0 & M_{zz_t} &= V_{zz_t} \\ 
    M_{xx_t} &= Q_{xx_t} + F_{x_t}^\intercal V_{xx_{t+1}} F_{x_t} &  M_{uu_t} &= Q_{uu_t} + F_{u_t}^\intercal V_{xx_{t+1}} F_{u_t} \\
    M_{ux_t} &= Q_{ux_t} + F_{u_t}^\intercal V_{xx_{t+1}} F_{x_t}  &  N_{x_t} &=  H_{x_{t+1}}F_{x_t} \\ 
    N_{u_t} &= H_{x_{t+1}} F_{u_t} & n_{1_t} &= H_{x_{t+1}}f_{1_t} + h_{1_{t+1}} \\
    \ N_{z_t} &= H_{z_t} & &  \
\end{align*}

Here, $Z_{w_t}$ is chosen such that the columns form an ortho-normal basis for the null-space of $N_{u_t}$, and $P_{y_t}$ is chosen such that its columns form an ortho-normal basis for the range-space of $N_{u_t}^\intercal$.  Hence $N_{u_t}$ and $P_{y_t}$ are orthogonal and their columns together span $\mathbb{R}^m$. 

Problem (\ref{opt:infeasible}), which is an unconstrained problem, has a solution given by the following: 
\begin{align}
	y_t^* &= -(N_{u_t}P_{y_t})^{\dagger} (N_{x_t} x_t + N_{z_t} x_{term} + n_{1_t}) \label{eq:controlvupdate} \\
	w_t^* &= -(Z_{w_t}^\intercal M_{uu_t} Z_{w_t})^{-1} Z_{w_t}^\intercal (M_{ux_t} x_t + M_{zu_t}^\intercal x_{term} + m_{u1_t}). \label{eq:controlwupdate} 
\end{align}
 
 Here $\dagger$ represents the pseudo-inverse. In the case that $P_{y_t}$ is a size-zero matrix (i.e. $\text{dim}(\text{null}(N_{u_t}))=m$), $Z_{w_t} = I_m$ (Identity matrix $\in \mathbb{R}^{m\times m}$), and $y_t$ has dimension 0. Correspondingly, when $Z_{w_t}$ is a size-zero matrix, $P_{y_t} = I_m$ and $w_t$ has dimension 0. Therefore we ignore one of the correspondingly empty update equations above (\ref{eq:controlvupdate} or \ref{eq:controlwupdate}). 
 
The control $u_t$ can now be expressed in closed-form as an affine function of the state $x_t$:
 \begin{align}
 	u_t^* &= K_{x_t} x_t + K_{z_t} x_{term} + k_{1_t}  \label{eq:controlpolicy} \\
	K_{x_t} &= -( P_{y_t}(N_{u_t}P_{y_t})^\dagger N_{x_t} + Z_{w_t} (Z_{w_t}^\intercal M_{uu_t} Z_{w_t})^{-1}Z_{w_t}^\intercal M_{ux_t} ) \label{eq:Kx} \\
	K_{z_t} &= -( P_{y_t}(N_{u_t}P_{y_t})^\dagger N_{z_t} + Z_{w_t} (Z_{w_t}^\intercal M_{uu_t} Z_{w_t})^{-1}Z_{w_t}^\intercal M_{zu_t} ) \label{eq:Kz} \\
	k_{1_t} &= -( P_{y_t}(N_{u_t}P_{y_t})^\dagger n_{1_t} + Z_{w_t} (Z_{w_t}^\intercal M_{uu_t} Z_{w_t})^{-1}Z_{w_t}^\intercal m_{u1_t} ) \label{eq:k1}
 \end{align}
Since the control is a function of the state, we can also express the constraint residual, i.e. the value of the term (\ref{eq:leastsquareswords}), as a function of the state. We can update the constraint-to-go to be this constraint residual, requiring preceding controls to enforce the residual is zero. We substitute  (\ref{eq:controlpolicy}-\ref{eq:k1}) into the right-hand side of  (\ref{eq:leastsquareswords}) to obtain:
\begin{align}
\begin{aligned}
	constraint\_to\_go_t(x_t) &= N_{x_t}x_t + N_{z_t}x_{term} + n_{1_t} -  \\
	 \ \ \ N_{u_t}P_{y_t}(N_{u_t}P_{y_t})^\dagger&(N_{x_t}x_t + N_{z_t}x_{term} + n_{1_t}) 
\end{aligned}
\end{align}
This results in the update for the terms $H_{x_t}$, $H_{z_t}$ and $h_{1_t}$:
\begin{align}
H_{x_t} &= (I - N_{u_t}P_{y_t}(N_{u_t}P_{y_t})^\dagger)N_{x_t} \\
H_{z_t} &= (I - N_{u_t}P_{y_t}(N_{u_t}P_{y_t})^\dagger)N_{z_t} \\
h_{1_t} &= (I - N_{u_t}P_{y_t}(N_{u_t}P_{y_t})^\dagger)n_{1_t}.
\end{align}	

 
The expression for the control is plugged in to the objective function of our optimization problem (\ref{opt:infeasible}) to obtain an update on the cost-to-go as a function of the state and $x_{init}$:
\begin{align}
    V_{xx_t} &=  M_{xx_t} + 2 M_{ux_t}^\intercal K_{x_t} + K_{x_t}^\intercal M_{uu_t} K_{x_t}  \label{update:Vxx} \\
    V_{zz_t} &= M_{zz_t} +  2 M_{ux_t}^\intercal K_{z_t} + K_{z_t}^\intercal M_{uu_t} K_{z_t} \label{update:Vzz} \\
    V_{zx_t} &= M_{zx_t} + M_{zu_t} K_{x_t} + K_{z_t}^\intercal M_{ux_t} + K_{z_t}^\intercal M_{uu_t} K_{x_t} \\
    v_{x1_t} &= m_{x1_t} + K_{x_t}^\intercal m_{u1_t} + (M_{ux_t}^\intercal + K_{x_t}^\intercal M_{uu_t}) k_{1_t} \label{update:vx} \\
    v_{z1_t} &= m_{z1_t} + M_{zu_t} k_{1_t} + K_{z_t}^\intercal m_{u1_t} + K_{z_t}^\intercal M_{uu_t} k_{1_t}. 
\end{align}

We have now presented forms for all terms in the $cost\_to\_go_t$ and $constraint\_to\_go_t$ functions, and computed linear feedback terms $K_{x_t}$, $K_{z_t}$ and $k_{1_t}$ in the process. This procedure can be repeated for time all time back to $t=0$. If $H_{x_0}x_{init} + H_{z_0}x_{term}+h_{1_0}=0$, then the sequence of control policies $\{K_{x_t}, K_{z_t}, k_{1_t}\}_{t\in\{0,...,T-1\}}$ will produce a sequence of states and controls that are optimal for our original problem (\ref{opt:globalwords}). Otherwise, no feasible solution exists. 

\subsubsection*{Trajectory}

Once we have computed all of the control policies for the entirety of the trajectory, the system can be forward simulated from the $x_{init}$ to obtain a form for all states and controls of the trajectory as functions of the variables $x_{init}$ and $x_{term}$.  We have:
\begin{align}
	x_t^* &= R_{a_t} x_{init} + R_{z_t} x_{term} + r_{1_t} \label{eq:xopt}\\
	u_t^* &= S_{a_t} x_{init} + S_{z_t} x_{term} + s_{1_t} \label{eq:uopt}
\end{align}
Here  the terms are initilaized as $R_{a_0} = I$,  $R_{z_0} = 0$, and $r_{1_t} = 0$. Otherwise, terms are recursively defined as:
\begin{align*}
	S_{a_t} &= K_{x_t} R_{a_t} & S_{z_t} &= K_{x_t} R_{z_t} + K_{z_t} \\
	 s_{1_t} &= K_{x_t} r_{1_t}  + k_{1_t} &  R_{a_{t+1}} &= F_{x_t} R_{a_t} + F_{u_t} S_{a_t} \\
	 R_{z_{t+1}} &= F_{x_t} R_{z_t} + F_{u_t} S_{z_t} & r_{1_t} &= F_{x_t} r_{1_t} + F_{u_t} s_{1_t} + f_{1_t}
\end{align*}

\subsubsection*{Lagrange Multipliers} 
Similar to the trajectory, the Lagrange multipliers associated with constraints (\ref{eq:dynamicswords3}-\ref{eq:term3}) as functions of $x_{init}$ and $x_{term}$. Writing the KKT conditions of the problem (\ref{opt:globalwords3}) that involve the multipliers, we have:

\begin{align}
&\begin{aligned}
\begin{bmatrix} Q_{xx_t} & Q_{ux_t}^\intercal \\  Q_{ux_t} & Q_{uu_t} \end{bmatrix} \begin{bmatrix} x_t^* \\ u_t^* \end{bmatrix} + \begin{bmatrix} q_{x1_t} \\ q_{u1_t} \end{bmatrix} + \begin{bmatrix} I & -F_{x_t}^\intercal \\ & -F_{u_t}^\intercal \end{bmatrix} \begin{bmatrix} \lambda_t \\ \lambda_{t+1} \end{bmatrix} = \begin{bmatrix} 0 \\ 0 \end{bmatrix} \\
 \ \ \ \ \ \ \ \ \  \forall t \in \{0,...,T-1\} 
\end{aligned} \label{eq:mults} \\
&Q_{xx_T} x_T^* + q_{x1_T} + \lambda_{T} + \mu_T = 0 \ \ \ \ \ \ \ \ \ \ \ \ \ \ \ \ \ \ \ \ \ \ \ \ \ \label{eq:mu}
\end{align}

As before, $\lambda_0 \in \mathbb{R}^n$ is the Lagrange multiplier associated with the constraint (\ref{eq:init3}), $\lambda_{t+1} \in \mathbb{R}^n, 0 \leq t < T$ are the Lagrange multipliers associated with constraints (\ref{eq:dynamicswords3}), and $\mu_T \in \mathbb{R}^n$ is the Lagrange multiplier associated with constraint (\ref{eq:term3}). 

Since the values of all optimal state and control variables, are known from (\ref{eq:xopt},\ref{eq:uopt}), the resulting system has $(T-1)(n+m) + n$ equations with $(T+1)n$ unknowns, and is an over-determined system. Stacking these conditions into a large matrix equation gives:  
\begin{align}
 A^\intercal \lambda = b \label{eq:multsglobal}
\end{align}
Here $A$ places all coefficients multiplying the Lagrange multipliers in the equations (\ref{eq:mults} and \ref{eq:mu}),  $b$ stacks all terms which do not depend on the multipliers, and $\lambda$ stacks all multipliers (including $\mu_T$). Since problem (\ref{opt:globalwords3}) is assumed to be convex and that the constraints are all linearly independent, there exists a unique solution in the dual space of the problem (\ref{opt:globalwords3}) \cite{boyd2004convex}. Furthermore, the matrix $A^\intercal$ is full column rank. This means that any combination of multipliers such that (\ref{eq:multsglobal}) is satisfied will correspond to the unique optimal set of multipliers $\lambda$. 

Therefore finding the least-squares solution of (\ref{eq:multsglobal}) will result in the multipliers per the discussion above. We have then that 
\begin{align}
AA^\intercal \lambda &= Ab  \label{eq:pseudoinverse}.
\end{align}
Looking at terms given in (\ref{eq:pseudoinverse}), we have: 
\begin{align}
	AA^\intercal &= \begin{bmatrix} 
		I & -F_{x_0}^\intercal \\
		-F_{x_0} & \Sigma_0 & -F_{x_1}^\intercal \\
		 & \ddots & \ddots & \ddots \\
		 & & -F_{T-2} & \Sigma_{T-2} & -F_{T-1}^\intercal \\
		 & & & -F_{T-1} & \Sigma_{T-1} & I \\
		 & & & & I & I
	\end{bmatrix} \label{eq:bandedAA} \\
	Ab &= \begin{bmatrix} 
	D_{a_0} & D_{z_0} \\
	\vdots & \vdots \\
	D_{a_T} & D_{z_T} 
	\end{bmatrix} \begin{bmatrix} x_{init} \\ x_{term} \end{bmatrix} + \begin{bmatrix} d_{1_0} \\ \vdots \\ d_{1_T} \end{bmatrix}
\end{align}
Here we have made use of the terms: 
\begin{align}
	&\Sigma_t = I + F_{x_t}F_{x_t}^\intercal + F_{u_t}F_{u_t}^\intercal   \\
	 &D_{a_0} = -Q_{xx_0}R_{a_0} - Q_{ux_0}^\intercal S_{a_0}  \\
	  &D_{z_0} = -Q_{xx_0}R_{z_0} - Q_{ux_0}^\intercal S_{z_0}   \\
	   &d_{1_0} = -Q_{xx_0}r_{1_0} - Q_{ux_0}^\intercal s_{1_0} - q_{x1_0}  \\
	   &\begin{aligned}
	  D_{a_{t+1}} = -Q_{xx_{t+1}}R_{a_{t+1}} - Q_{ux_{t+1}}^\intercal S_{a_{t+1}} + \ \ \ \ \ \ \ \ \ \ \  \ \ \ \  \ \  \\
	  \ \ \ \ \ \ \  F_{x_t}(Q_{xx_t}R_{a_t} + Q_{ux_t}^\intercal S_{a_t}) +
	  F_{u_t}(Q_{ux_t}R_{a_t} + Q_{uu_t}S_{a_t}) \end{aligned} \\
	  &\begin{aligned}
	  D_{z_{t+1}} = -Q_{xx_{t+1}}R_{z_{t+1}} - Q_{ux_{t+1}}^\intercal S_{z_{t+1}} + \ \ \ \ \ \ \ \ \ \ \  \ \ \ \ \ \ \  \\
\ \ \ \ \ \ \  F_{x_t}(Q_{xx_t}R_{z_t} + Q_{ux_t}^\intercal S_{z_t}) + 
	  F_{u_t}(Q_{ux_t}R_{z_t} + Q_{uu_t}S_{z_t}) \end{aligned} \\
	  &\begin{aligned}
	  d_{1_{t+1}} = -Q_{xx_{t+1}}r_{1_{t+1}} - Q_{ux_{t+1}}^\intercal s_{1_{t+1}} - q_{x1_{t+1}} + \\
	  F_{x_t}Q_{xx_t}r_{1_t} + F_{x_t}(Q_{ux_t}^\intercal s_{1_t} + q_{x1_t}) + \\
	  F_{u_t}(Q_{ux_t}r_{1_t} + Q_{uu_t}s_{1_t} + q_{u1_t}) \ \ \ \ \   \end{aligned} \\
	  &D_{a_{T}} = -Q_{xx_{T}}R_{a_{T}} \\
	  & D_{z_{T}} = -Q_{xx_{T}}R_{z_{T}} \\
	  & d_{1_{T}} = -Q_{xx_{T}}r_{1_{T}} - q_{x1_T}  
\end{align}
Notice that because (\ref{eq:bandedAA}) is a block-tridiagonal matrix, and the form of $Ab$ is explicit with respect to the variables $x_{init}$ and $x_{term}$, we can use a block substitution or sparse gaussian elimination method on the system (\ref{eq:pseudoinverse}) to determine the Lagrange multipliers as explicit functions of $x_{init}$ and $x_{term}$, in a computational complexity that scales linearly with respect to the length of the control horizon. 

Performing such a calculation will result in a simple expression for each of the multipliers of the form
\begin{align}
	\lambda_t &= L_{a_t}x_{init} + L_{z_t}x_{term} + l_{1_t}, \ \ t \in \{0,...,T\} \\
	\mu_T &= E_{a_T}x_{init} + E_{z_T}x_{term} + e_{1_T}.
\end{align}

Hence, all multipliers for the problem (\ref{opt:globalwords3}) were computed as explicit functions of the values $x_{init}$ and $x_{term}$.

\subsection{Method Recap}

In the previous two sections, the ingredients needed to parallelize Riccati recursion were developed. First, we showed that an LQR optimization can be broken into sub-problems, each of which correspond to end-point constrained LQR problems. A method was then presented for computing the solutions of the resulting sub-problems as functions of their end-points. Relating the Lagrange multipliers associated with the end-point constraints of the sub-systems resulted in a small system of equations that could be solved to determine the previously unknown link-points. Those solutions can then be redistributed back to the sub-problems so that the numerical values of their solutions can be computed. 

\section{EVALUATION} \label{sec:compare}
In this section we evaluate the computation time and analyze the feedback policies generated using the method presented. In particular, we compare our method against a well known parallel banded-matrix solver applied to the KKT system of equations of the global problem. We also demonstrate that the method we presented indeed produces the optimal trajectory of the global problem. We discuss the property that even though the resulting optimal trajectories are equivalent, the feedback policies generated by our method are slightly different. We illustrate this relationship with a simple example.   

\subsection{Parallel vs. Serial Riccati Recursion}
\begin{figure}[htb] 
  \centering
  \includegraphics[scale=0.17]{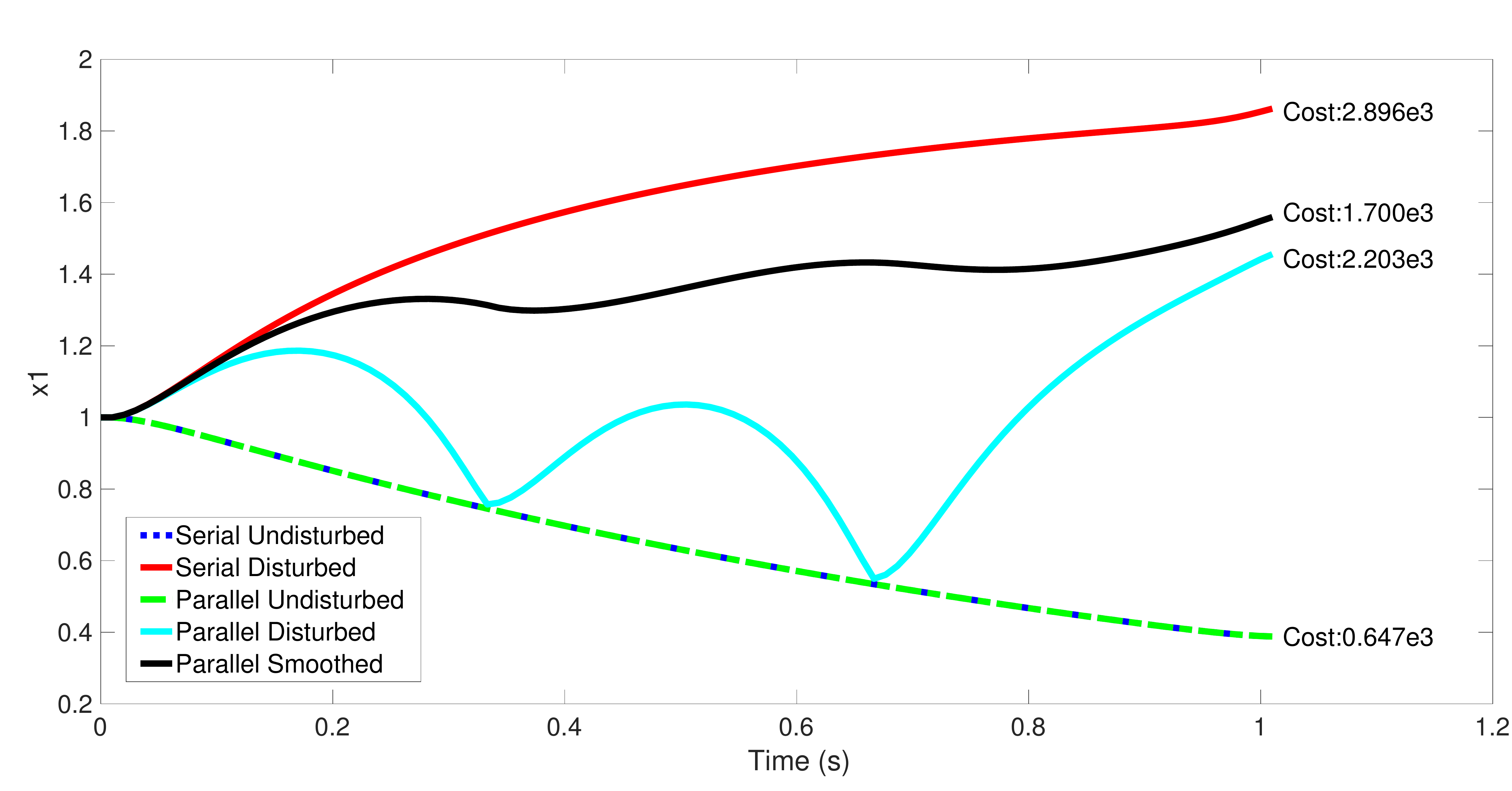} 
  \caption{Comparison of Serial and Parallel LQR Controllers from Problem (\ref{opt:doubleint}) on the disturbed system (\ref{eq:disturbed}).}
  \label{figure:parallelcompare} 
\end{figure}

Consider the following LQR problem: 
\begin{subequations}
\begin{align}
	\min_{u_0,...,u_{T-1}}& \  \alpha_T\|x_T\|_2^2 + \sum_{t=0}^{T-1} \alpha_t\|x_t\|_2^2 + \beta_t \|u_t\|_2^2  \\
	\text{s.t.} \ \ \ x_{t+1} &= \begin{bmatrix} 1 & dt \\ 0 & 1 \end{bmatrix} x_t + \begin{bmatrix} 0 \\ dt \end{bmatrix} u_t  \\
	x_0 &= \begin{bmatrix} 1& 0\end{bmatrix} ^\intercal 
\end{align} \label{opt:doubleint}
\end{subequations}
Here $\alpha_t=10$, $\alpha_T=10^3$, and $\beta_t=10^{-2}$. We compute the LQR feedback policies using standard Riccati recursion, as well as using the parallel method outlined above with $J=3$ sub-trajectories. Figure 2 demonstrates that, as we expect, the resulting optimal trajectories for both methods are identical (labels \say{Parallel Undisturbed} and \say{Serial Undisturbed}). We also apply the feedback policies generated in both methods to a disturbed system given by the following dynamics :
\begin{equation}
	x_{t+1} = \begin{bmatrix} 1 & dt \\ 0 & 1 \end{bmatrix} x_t + \begin{bmatrix} 0 \\ dt \end{bmatrix} u_t + \begin{bmatrix} 0 \\ x_t(1) / (x_t(1)^2 + 10^{-4}) \end{bmatrix} \label{eq:disturbed}
\end{equation}
We see in Figure 2 that when applied to the disturbed systems, the feedback policies now generate different solutions. Because each sub-trajectory of the parallel version are constrained to end at the \textit{would-be-optimal} point, we see strong corrections to these points. When the cost of these two trajectories are computed, we have $cost_s = 2.896e3$, and $cost_p=2.203e3$ for the serial and parallel policies, respectively. This behavior is interesting. On one hand, such correction can induce added stability to the policy, and result in improved performance as seen in this toy example. On the other hand, this may not always be beneficial. An optimal point for the assumed system dynamics may not be feasible for the true system, and as such a constraint forcing to one of these points might force the controller to apply unnecessarily strong actions in attempting to achieve the infeasible constraint. 

A middle-ground between these two methods can be obtained by re-computing the feedback policies for the first $J-1$ subtrajectories using the value function $cost\_to\_go_{\hat{t}_{j+1}}$ (value function at beginning of subtrajectory $j+1$) as an initialization for the $cost\_to\_go$ function at the end of subtrajectory $j$. In the second pass, the sub-trajectories no longer need to be constrained to terminate at the link-points. This has the effect of smoothing the control policies, and the result can also be seen in Figure \ref{figure:parallelcompare}. Note that this smoothing still results in the same optimal trajectory as the serial and parallel case when executed on the undisturbed dynamics. The cost of this smoothing procedure is roughly twice that of the non-smoothed parallel method, but might still be much faster than the serial method in some cases. 

\subsection{Computation Time vs. ScaLAPACK}
\begin{table}[htb]
\begin{tabular}{|c|c|c|c|c|c|}
\hline
n  & m  & T    & Cores & ScaLAPACK (s) & Parallel Riccati (s) \\ \hline
40 & 10 & 1024 & 32    & 0.074         & 0.043                \\ \hline
40 & 10 & 1024 & 64    & 0.057         & 0.037                \\ \hline
40 & 10 & 1024 & 1      & 0.797         & 0.727                \\ \hline
40 & 10 & 2048 & 32    & 0.116         & 0.074                \\ \hline
40 & 10 & 2048 & 64    & 0.085         & 0.057                \\ \hline
\end{tabular}
\caption{Comparing our method against ScaLAPACK's parallel banded matrix solver}
\label{table:compare}
\end{table}

We compare the computation time for our parallel method against a parallel solver for banded matrices. The KKT system of equations for an LQR problem form a banded coefficient matrix \cite{wright1996applying} which can be solved using specialized solvers for such systems. ScaLAPACK \cite{blackford1997scalapack} is a package which provides fortran routines for solving linear algebra problems on distributed processors.  It provides a method 'PDGBSV' which solves general banded systems using distributed processors and memory. Taking an approach based on the parallelization of banded systems was proposed in \cite{wright1990solution}. 

Table \ref{table:compare} compares the computation time of our method against the ScaLAPACK method. Here the dimensions of the state and control match that of the Cassie Robot, built by Agility Robotics. (Recall $n$ is the dimension of the state, and $m$ is the dimension of the control). All times are taken as the minimum over 10 runs on a Xeon Phi 1.3GHz 64-core processor. We use very long trajectory lengths to demonstrate the true effectiveness of parallelized methods. Our method outperforms the routine provided by ScaLAPACK in every case for this problem.  

\subsection{Other Methods}
The authors in \cite{farshidian2017real} present a method which also parallelizes the computation of solutions to LQR problems. Their method is based on using approximations of the value function, grabbed from previous solutions of the control problem when used in a Model Predictive Control framework, to avoid performing the dynamic programming on the entire length of the trajectory. One issue with such a method is that the resulting solutions are only approximate, and information will take at least $J$ cycles of planning to propagate to the start of the trajectory, where $J$ is the number of sub-trajectories. 

There exist other methods as well for parallelizing computation of controllers for LQR problems, such as in \cite{benner2008solving}. However, the approach in that work focuses on time-invariant, infinite-horizon problems and is therefore is not applicable for use in trajectory optimization for nonlinear systems with time-varying objectives. 

\section{Conclusion}
We have developed a method for parallelizing Riccati recursion, allowing for very efficient computation of feedback control policies for LQR problems. An interesting relationship between Lagrange multipliers of neighboring trajectories was developed, allowing for the distribution of optimization to multiple processors. We used a simple example to illustrate the differences between the serial and parallel versions of Riccati recursion, and presented some examples showing computational speedup using our method compared to existing methods. 
 









\bibliographystyle{./IEEEtran} 
\bibliography{root.bib}

\end{document}